\documentclass[11pt]{article}
\usepackage{amssymb}

\begin{document}

\newcommand{\e}{\varepsilon}
\renewcommand{\a}{\alpha}
\renewcommand{\b}{\beta}
\newcommand{\om}{\omega}
\newcommand{\Om}{\Omega}
\newcommand{\p}{\partial}
\renewcommand{\phi}{\varphi}

\newcommand{\N}{{\mathbb N}}
\newcommand{\R}{{\mathbb R}}
\newcommand{\EX}{{\mathbb E }}
\newcommand{\PX}{{\mathbb P }}

\newcommand{\cF}{{\cal F}}
\newcommand{\cG}{{\cal G}}

\newtheorem{theorem}{Theorem}
\newtheorem{lemma}{Lemma}
\newtheorem{remark}{Remark}

\title{Enstrophy Dynamics of Stochastically Forced Large-Scale
Geophysical Flows}

{\tiny
\author{Dirk Bl\"omker$^1$, Jinqiao Duan$^2$,
and Thomas Wanner$^3$\\
\\
1.  Institut f\"ur Mathematik \\
    Universit\"at Augsburg\\
    86135 Augsburg, Germany \\
2.  Department of Applied Mathematics \\
Illinois Institute of Technology \\
Chicago, IL 60616, USA  \\
3.  Department of Mathematics and Statistics\\
   University of Maryland, Baltimore County\\
   Baltimore, MD 21250, USA}
}

\date{June 28, 2001}

\maketitle

\begin{abstract}

Enstrophy is an averaged measure of fluid vorticity. This quantity is
particularly important in {\em rotating} geophysical flows. We
investigate the dynamical  evolution  of enstrophy for large-scale quasi-geostrophic
flows under random wind forcing. We obtain upper bounds on the
enstrophy, as well as results establishing its H\"older continuity and
describing the small-time asymptotics.

\medskip
{\bf Key Words:} Enstrophy, geophysical fluid dynamics, stochastic
partial differential equations

\medskip
{\bf Mathematics Subject Classifications (2000)}:
37A25, 60H15, 76D05, 86A05

\end{abstract}

\newpage

\section{Introduction}

Randomness is ubiquitous in fluid systems.  Macroscopic partial
differential equation models for fluid flows contain such randomness
as stochastic forcing, uncertain parameters, random sources, and random
boundary conditions.

There has been active recent research on stochastic approaches to
geophysical flows \cite{Muller, Holloway,Griffa, Samelson,
Br-Du-Wa:00} and numerical simulations of stochastically forced
geophysical flows \cite{Sura, Sura2, Brannan1, Brannan2, SD}.  It is
generally understood that random fluctuations can have delicate impact
on geophysical fluid dynamics \cite{Muller, Griffa, Sura, Sura2,
DuanSchm}.

A class of large-scale geophysical flows under random forcing are
modeled by the quasi-geostrophic equation \cite{Muller}:
\begin{equation}
   \Delta \psi_t + J(\psi, \Delta \psi ) + \beta \psi_x
   =\nu \Delta^2 \psi - r \Delta \psi + \dot{W}  \; ,
   \label{qg}
\end{equation}
where $\psi(x,y,t)$ is the stream function ($\psi_x:=\partial_x\psi$),
$\beta \geq 0$ is the meridional gradient of the Coriolis parameter,
$\nu >0$ is the viscous dissipation constant, $r>0$ is the Ekman
dissipation constant, and $W(x,y,t)$ is a space-time Wiener process to
be defined below on a probability space $(\Omega, {\cal A},
\PX)$. Moreover, $J(f, g)=f_x g_y -f_y g_x$ denotes the Jacobian
operator. The generalized time derivative $\dot{W}$ models the noisy
wind forcing.

Introducing $\om (x,y,t) = \Delta \psi(x, y,t)$, equation (\ref{qg})
can be rewritten in the form
\begin{equation}
       \om_t + J(\psi, \om ) + \beta \psi_x
        =\nu \Delta \om - r \om  +  \dot{W} \; ,
 \label{eqn}
\end{equation}
where $(x, y) \in D$ and $D \subset \R^2$ denotes a bounded domain
with sufficiently smooth boundary. The boundary conditions are no
normal flow ($\psi=0$ on $\partial D$) and free-slip ($\om=0$ on
$\partial D$) as in Pedlosky \cite[p.~34]{Ped96} or in Dymnikov and
Kazantsev \cite{Dymnikov}:
\begin{eqnarray}
 \psi = \om =  0  \quad \mbox{on} \; \partial D \; . \label{BC}
\end{eqnarray}
An appropriate initial condition  $\om(0)$ is also imposed.

The mean enstrophy for a fluid flow is half the squared mean-square
norm of the vorticity \cite{Pedlosky, Sal98}, i.e., we have
$$
Ens(t) = \frac12 \cdot \EX \int_D |\om(x,y,t)|^2 d(x,y) \; .
$$
The enstrophy $Ens(t)$ is an averaged measure of fluid vorticity
$\om(t)$. In this paper, we discuss the time evolution of the
enstrophy. We present results which establish upper bounds on the
enstrophy, as well as results on H\"older continuity and small-time
asymptotics for $Ens(t)$. These results are contained in Sections 3,
4, and 5, respectively. The mathematical framework for our discussion
is described in Section 2.


\section{Mathematical Framework}

As it stands, the stochastic quasi-geostrophic equation (\ref{eqn})
still has to be given a mathematically precise meaning. This can be
accomplished using the framework of stochastic partial differential
equations, as described for example in \cite{DaPrato}.  In our
situation, we formally rewrite (\ref{eqn}) in the Ito formulation
\begin{equation}
    d\om =(\nu \Delta \om - r \om -\beta \psi_x -J(\psi, \om ))dt + dW \; .
    \label{neweqn}
\end{equation}
In the following we use the abbreviations $L^2=L^2(D)$,
$L^{\infty}=L^{\infty}(D)$, $H^k_0=H^k_0(D)$, $H^k =H^k(D)$, $ 0 < k <
\infty$, for the standard Sobolev spaces. Let $<\cdot, \cdot>$ and $\|
\cdot \| $ denote the standard scalar product and norm in $L^2$,
respectively. Moreover, the norms for $H^k_0$ and $L^{\infty}$ are
denoted by $\| \cdot \|_{H^k}$ and $\|\cdot\|_{\infty}$,
respectively. Due to the Poincar\'e inequality
\cite[p.~164]{Gilbarg-Trudinger}, the expression $\| \Delta \cdot \|$
is an equivalent norm for $H^2_0$.  It is well-known that the operator
$A = \nu \Delta : L^2 \to L^2$ with domain $D(A) = H^2 \cap H_0^1$ is
self-adjoint.  Note that $A$ generates an analytic semigroup $S(t)$ on
$L^2$ (\cite{Pazy}).  The spectrum of $ A$ consists of eigenvalues $0
> \lambda_1 > \lambda_2 \ge \lambda_3 \ge \ldots$ with corresponding
normalized eigenfunctions $\phi_1$, $\phi_2$, $\ldots$. The set of
these eigenfunctions is complete in $L^2$. For example, for the square
domain $D = (0,1) \times (0,1)$ the eigenvalues are given by
$-\nu(m^2+n^2)\pi^2$ for $m,n \in \N$, and the associated
eigenfunctions are suitable multiples of $\sin(m\pi x) \sin(n\pi y)$.

Now we can define an appropriate class of Wiener processes $W$. Let
$\beta_k$, $k \in \N$, denote a family of independent real-valued
standard Brownian motions. Furthermore, choose positive constants
$\mu_k$, $k
\in \N$, such that
\begin{equation}  \label{As-Wsum}
  \sum_{k=1}^\infty \frac{\mu_k^2}{|\lambda_k|^{1-\theta}} < \infty
\end{equation}
for some $0 < \theta < 1$. Then we consider the Wiener process $W$ defined
by
\begin{eqnarray}
  W(t) := \sum_{k=1}^\infty  \mu_k  \cdot \beta_k(t) 
          \cdot \phi_k \; ,
  \quad t \ge 0 \; .
  \label{Wiener}
\end{eqnarray}
Note that we explicitly allow Wiener processes $W$ whose covariance
operator is not of trace-class, i.e., for which $\sum_{k=1}^\infty
\mu_k^2=\infty$.


For the domain $D$ we basically assume that the eigenfunctions
$\phi_k$ of $A$ satisfy
\begin{equation} \label{bound-EF}
  \phi_k \in C_0(\bar{D}), \quad |\phi_k(x,y)| \leq C,
\end{equation}

\begin{displaymath}
  |\nabla \phi_k(x,y)| \leq C \sqrt{|\lambda_k|},
\end{displaymath}
for all $(x,y) \in D$ and $k \in \N$, where $C>0$ denotes a constant
which depends only on~$D$. Domains~$D$ which satisfy these conditions
include rectangular domains, as well as equilateral triangles.
Unfortunately, there are many domains for which they are violated. See
for example \cite{aurich:etal:99a}. However, in this paper it is
conjectured that in (\ref{bound-EF}) one generally should expect an
upper bound which is logarithmic in $|\lambda_k|$. Even though our
results remain valid in this situation, we will assume the above
stronger condition.

Under the above assumptions, Theorem 5.2.9 in \cite{DaPrato2}
guarantees that the stochastic convolution
\begin{eqnarray}
   W_A(t)  =  \int_0^t S(t-s)dW(s) \; , \quad t>0 \; ,
  \label{conv}
\end{eqnarray}
has a continuous version with values in $C_0(D)$, the Banach space of
continuous functions satisfying zero Dirichlet boundary conditions on
$D$. To be more precise, $W_A$ has a version which is even H\"older
continuous with some small exponent, which depends on the asymptotic
behavior of the coefficients~$\mu_k$.

If we define the nonlinear operator $F$ by $F(\om) = -r \om
-\beta\psi_x -J(\psi, \om )$, then (\ref{neweqn}) can be rewritten as
the abstract evolution equation together with initial condition
\begin{eqnarray}
    d\om & = & (A \om + F(\om) )dt + dW \; ,\label{QG}
     \\
    & & \om (0) = \om_0.\nonumber
\end{eqnarray}
For technical reasons we translate the operator~$A$ by a suitable
multiple of the identity. Consider a constant $\alpha \ge 0$ which
will be chosen later on, usually sufficiently large.  Defining
$A_\alpha := A - \alpha I$, we get the initial value problem
\begin{eqnarray}
    d\om & = & (A_\alpha \om + F(\om)+\alpha\om )dt + dW \; ,
      \nonumber
      \\
    & & \om (0) = \om_0
        \nonumber
\end{eqnarray}
or in mild (integral) form
\begin{eqnarray}
     \om(t) = S_\alpha(t) \om_0
     + \int_0^t S_\alpha(t-s) \left( F(\om (s)  )+\alpha\om(s)
        \right) ds
     + W_{A_\alpha}(t) \; ,
      \label{mild}
\end{eqnarray}
where the analytic semigroup $S_\alpha$ is given by $S_\alpha(t) =
e^{-t\alpha} \cdot S(t)$ for $t>0$ and the stochastic convolution
$W_{A_\alpha}(t)$ is defined as in (\ref{conv}) with the semigroup~$S$
replaced by~$S_\alpha$. Finally, let $U := \om -
W_{A_\alpha}$. Then~$U$ is the weak solution of
\begin{eqnarray}  \nonumber
    \partial_t U & = & A U + F(U+V)+\alpha V \; , \\
     & & U (0) = \om_0,  \nonumber
\end{eqnarray}
where we use the abbreviation $V := W_{A_\alpha}$. Notice that both
$U$ and $V$ depend on $\alpha$.


\section{Enstrophy Estimate: Upper Bounds}

We begin by establishing upper bounds on the time evolution of the
enstrophy $Ens(t) = \EX \|\om(t)\|^2 / 2$.  Improving the a-priori
estimate of \cite[Section~3]{Br-Du-Wa:00} we obtain
\begin{eqnarray}
 \lefteqn{ \frac12 \cdot \frac{d}{dt}\|U(t)\|^2
 \leq   \|\nabla U(t)\|^2\cdot(\e-\nu)}
\nonumber\\&& \displaystyle
  + \; \|U(t)\|^2\cdot \left( (\e-r+c_1\b) + 
  C \cdot \|V\|_{\infty} \cdot \left( 1+C_\e \cdot \|V\|_{\infty}
    \right) \right) \nonumber\\[1ex]&& \displaystyle
  + \; C_\e \cdot \left( r + \b + \alpha^2 \right) \cdot
      \|V\|_{\infty}^2 +
    C \cdot \|V\|_{\infty}^3 + C_\e \cdot \|V\|_{\infty}^4,
\end{eqnarray}
where $C$ denotes a generic constant which depends only on $D$, and
whose specific value may change from line to line. Similarly, $C_\e$
denotes a generic constant which depends only on $D$ and $\e$, where
$\e>0$ is some arbitrarily small number. The constant $c_1$ denotes
the optimal constant in the Poincare inequality $\| U \| \le c_1 \|
\nabla U \|$ for mean zero functions~$U$. For $\e<\nu$ the improved
a-priori estimate immediately yields the following lemma.

\begin{lemma}  \label{lem-APR}
For any $\gamma > - \nu \cdot c_1^{-2} - r + c_1 \cdot \b$ there exist
constants depending only on $\gamma$, $D$, $\b$, and $r$ which are all
denoted by $C$ such that
\begin{equation}
\frac{d}{dt}\|U(t)\|^2
 \leq A(t) \|U(t)\|^2+B(t), \; \quad t\ge0
\end{equation}
with
\begin{eqnarray}  \label{Def-AB}
  A(t) & = & 2 \gamma + C \cdot \left( \|V\|_{\infty} +
    \|V\|_{\infty}^2 \right), \nonumber \\
  B(t) & = & C \cdot \left( \left( 1+\alpha^2 \right) \cdot 
    \|V\|^2_{\infty} + \|V\|_{\infty}^4 \right) > 0.
\end{eqnarray}
Together with Theorem I.6.1 in \cite{Hale} this yields
\begin{equation}
  \|U(t)\|^2  \leq \|\om_0\|^2 \cdot e^{\int_0^t A(s) ds}
  +\int_0^t B(s) e^{\int_s^t A(\tau)d\tau} ds,\; \quad t\ge0 .
  \label{bounded}
\end{equation}
\end{lemma}

Suppose for simplicity that $\omega_0$ and $W$ are stochastically
independent. It is possible to drop this assumption in this
section, but in this case we additionally need
$\EX\|\om_0\|^{2+\delta}<\infty$ for some small $\delta>0$.

The critical term for taking the expectation in (\ref{bounded}) is the
squared $L^\infty$-norm of $V = W_{A_\alpha}$ in the exponent, which
is in general not finite. To complicate matters further, the
$L^\infty$-norm in the exponent cannot easily be dealt with, since we
do not have a Hilbert space structure.

For our situation we will improve on some ideas of
\cite{PhD-Gugg}. Using Fernique's Theorem \cite[Theorem~2.6]{DaPrato}
we get that $\PX(t\|V(\tau)\|^2_\infty>r^2)\le 1/(1+e^{1+32\lambda
r^2})$ implies $\EX (e^{\lambda t \| V(\tau)\|_{\infty}^2}) \le e^{16
\lambda r^2} + e^2 / (e^2-1)$ for any $t,\tau,r,\lambda>0$.  Hence, by
Jensen's inequality
$$
    \EX \left( e^{ \lambda \int_0^t\| V(\tau)\|_{\infty}^2 d\tau}
      \right)  \le
    \frac{1}{t} \cdot \int_0^t \EX \left( e^{ \lambda t \| 
      V(\tau)\|_{\infty}^2}\right) d\tau
    \leq C_\lambda,
$$
provided $\PX(t\|V(\tau)\|^2_\infty>1)\le 1/(1+e^{1+32\lambda})$ for
any $\tau\le t$. The latter inequality follows immediately from
Chebychev's inequality, provided we have
\begin{equation}  \label{BG1}
   t\cdot \EX \| V(\tau)\|_{\infty}^2\le \frac1{1+e^{1+32\lambda}}
\end{equation}
for any $\tau\le t$.  The following lemma is proven similarly as
\cite[Theorem~5.20]{DaPrato}.

\begin{lemma}  \label{supV}
For any $p\ge1$ and any sufficiently small $\theta>0$ there exists a
constant~$C$ which depends only on $p$, $\theta$, and $D$ such that
for any $\tau\ge0$
\begin{equation}  \label{V-bound}
  \EX \| V(\tau)\|_{\infty}^{2p} \le C\cdot
  \underbrace{\left(\sum_{k=1}^{\infty}\frac{\mu_k^2}{\alpha-\lambda_k}
  \cdot |\lambda_k|^{\theta}\right)^p}_{=:\phi(\alpha)^p}.
\end{equation}
\end{lemma}

We remark that the assumption (\ref{bound-EF}) on the eigenfunctions
is essential for the proof of this lemma. Notice also that the series
in (\ref{V-bound}) is finite according to
(\ref{As-Wsum}). Lemma~\ref{supV} implies that (\ref{BG1}) is
satisfied for any $t\le C/\phi(\alpha)$. It is now straightforward to
verify that
$$
\EX e^{m\int_s^tA(\tau)d\tau}\leq C \cdot e^{m(t-s)2\gamma}
$$
for any $s\le t\le C/\phi(\alpha)$ and $m=1,2$, where $\gamma$ and
$A(t)$ were defined in Lemma \ref{lem-APR}. Moreover, Lemma
\ref{supV} implies for $B(t)$ in (\ref{Def-AB})
$$
   \left( \EX B(t)^2 \right)^{1/2}
   \le C \cdot \left( \left( 1+\alpha^2 \right) \cdot
   \phi(\alpha) + \phi(\alpha)^2 \right)
$$
for any $t>0$. We finally obtain from (\ref{bounded}) that
\begin{eqnarray}
\EX \|U(t)\|^2
&\leq& \EX \left( \|\omega_0\|^2 \cdot e^{\int_0^tA(s)ds}
         +\int_0^tB(s)e^{\int_s^tA(\tau)d\tau}ds\right)\nonumber\\
&\leq& \EX \|\omega_0\|^2 \cdot \EX e^{\int_0^tA(s)ds}
         +\int_0^t \left( \EX B(s)^2 \right)^{1/2}
         \left( \EX e^{2\int_s^tA(\tau)d\tau} \right)^{1/2}
         ds\nonumber\\
&\leq& C \cdot \EX \|\omega_0\|^2 \cdot e^{2\gamma t} \nonumber\\
& & \quad + \, C \cdot \left( \left( 1+\alpha^2 \right) \cdot
         \phi(\alpha)+\phi(\alpha)^2 \right) \cdot
         \int_0^te^{2\gamma\tau}d\tau
\end{eqnarray}
for any $t\le C/\phi(\alpha)$. Using
$\EX\|\om\|^2\le2\EX\|U\|^2+2\EX\|V\|^2$, this immediately implies the
following theorem on upper bounds for the enstrophy.

\begin{theorem}[Upper Bound] \label{om-bound}
  Suppose that $\omega_0$ and $W$ are stochastically independent and
  that $\EX\|\omega_0\|^2<\infty$. Moreover, let $\gamma >
  -\nu \cdot c_1^{-2} - r + c_1 \cdot \b$.
  Then $Ens \in L^\infty([0,T])$ for any $T>0$. More precisely,
  \begin{eqnarray*}
    Ens(t) & \le &
    C \cdot \EX \|\omega_0\|^2 \cdot e^{2\gamma t} +
      C \cdot \varphi(\alpha) \\
    & & \,\, + \,\, C \cdot \left( \left( 1+\alpha^2 \right) \cdot
      \phi(\alpha) + \phi(\alpha)^2 \right) \cdot
      \int_0^te^{2\gamma\tau}d\tau
  \end{eqnarray*}
  for any $t\le C/\phi(\alpha)$, with constants $C$ independent of
  $t$, $\alpha$, and $\omega_0$.
\end{theorem}

\begin{remark} \label{om-bound-remark}
It can be shown that for any choice of $p\ge1$ similar bounds hold for
$\EX\sup_{\tau\in[0,t]}\|U(\tau)\|^{2p}$, provided both
$\EX\|\omega_0\|^{2p}<\infty$ and $t\le C/\phi(\alpha)$.

Moreover, all bounds on $\EX\|U(t)\|^{2p}$ or
$\EX\sup_{\tau\in[0,t]}\|U(\tau)\|^{2p}$ immediately imply analogous
bounds on $Ens(t)$ or $\EX\sup_{\tau\in[0,t]}\|\om(\tau)\|^{2p}$.  For
this one has to employ estimates for
$\EX\sup_{\tau\in[0,t]}\|W_{A_\alpha}(\tau)\|^{2p}$ which can be
obtained for example as in \cite[Corollary~2.3]{Gu-Ki-Ni:01}.
\end{remark}

In order to obtain a bound for $t \rightarrow \infty$, we note that
for $\alpha \rightarrow \infty$ one obviously has $\phi(\alpha)
\rightarrow 0$. However, the rate of convergence is essential. To this
end, we distinguish two cases. If we suppose that $\sum_{k=1}^\infty
\mu_k^2 \cdot |\lambda_k|^\theta < \infty$ for some $\theta > 0$, then
the estimate $\phi(\alpha) \le \sum_{k=1}^\infty \mu_k^2 \cdot
|\lambda_k|^\theta / \alpha$ is immediate --- and choosing~$\alpha$
proportional to~$t$ in Theorem~\ref{om-bound} furnishes
\begin{displaymath}
   Ens(t) \le  C \cdot \left( \EX \|\omega_0\|^2 \cdot e^{2\gamma t}
         + t \cdot \int_0^t e^{2\gamma\tau} d\tau + 1
   \right) \quad\mbox{ for all } t\ge0 \, .
\end{displaymath}
If on the other hand $\sum_{k=1}^\infty \mu_k^2 \cdot
|\lambda_k|^\theta = \infty$, we additionally assume $\mu_k^2 \le C
k^{-\mu}$ for some $\mu \in (\theta,1]$, with arbitrarily small
$\theta$ defined in Lemma \ref{supV}. (For $\mu > 1$ we can always
find some small $\theta$ such that the first case applies.) Using the
fact that $|\lambda_k| \sim Ck$ for $k\rightarrow\infty$
(cf. \cite{Evans}) we obtain
\begin{eqnarray*}
  \phi(\alpha) & \le & C \cdot \sum_{k=1}^\infty
    \frac{k^{-\mu}}{ck+\alpha} \cdot k^{\theta} \le
    C \cdot \int_0^\infty \frac{k^{\theta-\mu}}{ck+\alpha} \, dk \\
  & = & C \cdot \alpha^{\theta-\mu} \cdot
    \int_0^\infty \frac{\tau^{\theta-\mu}}{c\tau + 1} \, d\tau \; .
\end{eqnarray*}
Choosing $\alpha^{\mu-\theta}$ proportional to $t$ in Theorem
\ref{om-bound}, we derive
\begin{displaymath}
   Ens(t)\le  C \cdot \left( \EX \|\omega_0\|^2 \cdot e^{2\gamma t}
         + t^{2/(\mu-\theta)-1} \cdot \int_0^te^{2\gamma\tau} d\tau
         + 1 \right)
   \mbox{ for all } t\ge0 \; .
\end{displaymath}
Notice that $2/(\mu-\theta)-1 > 1$.  We have proved the following
result.

\begin{theorem}[Global Upper Bound] \label{thmglobal}
  Assume again that $\omega_0$ and $W$ are stochastically independent,
  that $\EX\|\omega_0\|^2 < \infty$, and let $\gamma >
  -\nu \cdot c_1^{-2} - r + c_1 \cdot \b$. Then the following
  holds.
  \begin{itemize}
  \item[(a)] If $\mu_k^2 \le C k^{-\mu}$ for some $\mu \in (0,1]$,
             we can find a $\tilde\mu \in (0,\mu)$ such that
\begin{displaymath}
   Ens(t) \le  C \cdot \left( \EX \|\omega_0\|^2 \cdot e^{2\gamma t}
     + t^{(2-\tilde\mu)/\tilde\mu} \cdot
     \int_0^t e^{2\gamma\tau} d\tau + 1 \right)
   \mbox{ for all } t \ge 0 \; .
\end{displaymath}
  \item[(b)] If $\sum_{k=1}^\infty \mu_k^2 \cdot |\lambda_k|^\theta
             <\infty$ for some $\theta>0$ then
\begin{displaymath}
   Ens(t) \le C \cdot \left( \EX \|\omega_0\|^2 \cdot e^{2\gamma t}
         + t \cdot \int_0^t e^{2\gamma\tau} d\tau + 1 \right)
   \mbox{ for all } t \ge 0 \; .
\end{displaymath}
  \end{itemize}
  The constants $C$ are independent of $t$ and $\om_0$, but they can
  depend on $\theta$, $\gamma$, $\mu$, $\tilde\mu$, the domain $D$,
  or the coefficients in (\ref{qg}).
\end{theorem}

\begin{remark}
   If $-\nu \cdot c_1^{-2} - r + c_1 \cdot \b \ge 0$, then
   necessarily $\gamma > 0$. In this case we obtain an exponentially
   growing upper bound for $Ens(t)$ with growth rate slightly larger
   than $-\nu \cdot c_1^{-2} - r + c_1 \cdot \b$.

   If on the other hand $-\nu \cdot c_1^{-2}-r+c_1 \cdot \b<0$, then we
   can choose $\gamma<0$. This furnishes a polynomial upper bound which
   grows at least linearly in time. The precise growth exponent is
   determined by the regularity of the noise.
\end{remark}

\begin{remark}
   As we already stated in the beginning of this section, one can
   remove the condition of stochastic independence of $\omega_0$
   and $W$ in the previous theorem, if one additionally assumes
   $\EX\|\omega_0\|^{2+\delta}<\infty$ for some small $\delta$.
\end{remark}

Our above results hold for a large class of noise processes, in
particular also for more irregular Wiener processes $W$ whose
covariance operator is not of trace-class.  If, however, one assumes
that the Wiener process is of trace-class, i.e., if $\mathrm{Tr}(Q) =
\sum_{k=1}^\infty \mu_k^2 < \infty$, then the results can be improved
significantly by employing Ito's formula. One advantage of this
approach is that it avoids the conditions on the eigenfunction in
(\ref{bound-EF}).  Therefore, we will briefly outline the main ideas.

By applying Ito's formula \cite[Section~4.5]{DaPrato} to the squared
$L^2$-norm of the vorticity $\omega(t)$, it can easily be verified
that
\begin{displaymath}
  \EX\|\om(t)\|^2 =
  2 \EX \int_0^t <A\om(\tau)+F(\om(\tau)),\om(\tau)> d\tau
  +\mathrm{Tr}(Q) \cdot t \, ,
\end{displaymath}
where $\mathrm{Tr}(Q)=\sum_{k=1}^\infty \mu_k^2$ denotes the trace of
the covariance operator~$Q$ of~$W$. Using calculations analogous to
the ones leading to the a-priori estimate in Lemma \ref{lem-APR}, we
formally obtain
\begin{eqnarray*}
  \partial_t \EX\|\om(t)\|^2 &=&
    2 \EX <A\om(\tau)+F(\om(\tau)),\om(\tau)> + \mathrm{Tr}(Q)\\
  &\le& 2\gamma \cdot \EX\|\om(t)\|^2 + \mathrm{Tr}(Q) \, ,
\end{eqnarray*}
where $\gamma > -\nu \cdot c_1^{-2} - r + c_1 \cdot \b$ as in
Lemma~\ref{lem-APR}. Hence, for any $t\ge0$ we have
\begin{equation} \label{trclest}
  Ens(t) \le Ens(0) \cdot e^{2\gamma \cdot t} +
  \mathrm{Tr}(Q) \cdot \frac{e^{2\gamma \cdot t}-1}{4\gamma}
  \; .
\end{equation}
Especially if one can choose a growth exponent $\gamma < 0$, this
significantly improves the estimates of Theorem~\ref{thmglobal}, since
in this case the right-hand side of~(\ref{trclest}) approaches
$-\mathrm{Tr}(Q) / (4 \gamma)$ for $t \to \infty$.

\section{Enstrophy Estimate: H\"older Continuity}

In this section we establish regularity properties of the enstrophy as
a function of time. More precisely, we will prove that $Ens(t) = \EX
\|\om(t)\|^2 / 2$ is H\"older continuous.  To this end, we need the
following lemma from \cite{Br-Du-Wa:00}.

\begin{lemma} \label{f_global}
  Define a nonlinear mapping $\cF : C( [0,T]; H_0^1) \to 
  C( [0,T]; L^2)$ by
  \begin{displaymath}
    (\cF(\omega))(t) := \int_0^t S(t-s) F(\omega(s)) ds \; ,
    \quad\mbox{ for }\quad t \in [0,T] \; ,
  \end{displaymath}
  where $\omega \in C( [0,T]; H_0^1)$, and $A$ and $F$ are as in
  (\ref{QG}). Then $\cF$ is continuous,
  and it can be extended to a continuous mapping from the space
  $C( [0,T]; L^2)$ to $C([0,T]; L^2)$. Furthermore, the image of
  the extended mapping $\cF$ is contained in $C([0,T],H^a(D))$
  for $0 \le a < 1/2$.
\end{lemma}

In fact, it is shown in \cite{Br-Du-Wa:00} that for arbitrary positive
constants $a \in [0,1/4)$ and $\rho \in (0,1/4)$ satisfying $0 < \rho
+ a < 1/4$ the estimate
\begin{eqnarray*}
\lefteqn{\left\| (-A)^a \int_0^t S(t-s) F(\om (s)) ds \right\|} \\
& \leq & \frac{rc + \beta c}{1-a} \cdot t^{1-a} \cdot
    \sup\limits_{0\leq s \leq t} \|\om(s)\|   \\
& & \quad + \left( \frac{8c}{1-4\rho-4a} \cdot t^{\frac14-\rho-a} +
    \frac{4c}{1-2\rho-2a} \cdot t^{\frac12-\rho-a} \right)
     \cdot \sup\limits_{0\leq s \leq t} \|\om(s)\|^2
\end{eqnarray*}
holds.  Especially for $a = 0$ we obtain
\begin{eqnarray*}
  \lefteqn{ \| \cF(\om)(t) \| =\left\| \int_0^t S(t-s) F(\om (s))
    ds \right\|} \\
  & \leq & (rc + \beta c) \cdot t \cdot
      \sup\limits_{0\leq s \leq t} \|\om(s)\|   \\
  & & \quad + \left( \frac{8c}{1-4\rho} \cdot t^{\frac14-\rho} +
      \frac{4c}{1-2\rho} \cdot t^{\frac12-\rho} \right)
      \cdot \sup\limits_{0\leq s \leq t} \|\om(s)\|^2,
\end{eqnarray*}
for every $0 < \rho < 1/4$. Together with Theorem \ref{om-bound} and
Remark \ref{om-bound-remark} these bounds immediately furnish the
following result.

\begin{lemma}  \label{F-bound}
  Suppose that $\omega_0$ and $W$ are stochastically independent and
  that for some $p\ge1$ we have $\EX\|\omega_0\|^{2p}<\infty$.
  Moreover, let $T>0$ be arbitrary, and
  let $a \in [0,1/4)$ and
  $\rho \in (0,1/4)$ be such that $0 < \rho + a < 1/4$. Then there
  exists a constant $C$ such that
  \begin{displaymath}
    \EX \| (-A)^a \cF(\om)(t) \|^p \le
    C \cdot t^{p \cdot \left( \frac{1}{4} - \rho - a \right)}
    \quad\mbox{ for all }\quad t \in [0,T] \; .
  \end{displaymath}
\end{lemma}

The following theorem states our main result on the regularity of the
enstrophy. It will be proved in the remainder of this section.

\begin{theorem}[H\"older Continuity] \label{holder}
  Suppose that $\omega_0$ and $W$ are stochastically independent and
  that $\EX\|\omega_0\|^{4}<\infty$. Then the enstrophy $Ens(\cdot)$
  defined by $Ens(t) = \EX \|\om(t)\|^2 / 2$ is H\"older continuous
  with arbitrary exponent less than $1/4$ on every compact interval
  in $(0,\infty)$.
\end{theorem}

\begin{remark}
  Note that in general we cannot expect the solution $\om$ to be
  H\"older continuous with arbitrary exponent less than $1/4$, since
  $W_A$ is in general less regular. As one can see from (\ref{mild})
  with $\alpha = 0$, one cannot expect $\om$ to be more regular than
  $W_A$.
\end{remark}

To prove the above theorem establishing the H\"older continuity of the
enstrophy, we first define
\begin{equation}\label{Def-G}
  \cG(\om)(t) := S(t)\om_0 + \cF(\om)(t).
\end{equation}
According to (\ref{mild}) for $\alpha = 0$, we therefore have $\om(t)
= \cG(\om)(t) + W_A(t)$. Consider a fixed interval $J =
[\varepsilon,T] \subset (0,\infty)$. For $t\in J$ and $h$ with $t+h\in
J$ the identity (\ref{mild}) then implies
\begin{eqnarray*}
  \lefteqn{ \EX \| \om(t+h)\|^2 - \EX\|\om(t)\|^2}\\
  &=& \EX\|\cG(\om)(t+h)+W_A(t+h)\|^2-\EX\|\cG(\om)(t)+W_A(t)\|^2\\
  &=& \EX\|\cG(\om)(t+h)\|^2-\EX\|\cG(\om)(t)\|^2\\
  && +2 \EX < \cG(\om)(t+h),W_A(t+h)>-2\EX<\cG(\om)(t),W_A(t)>\\
  && +\EX\|W_A(t+h)\|^2-\EX\|W_A(t)\|^2\\
  &=:& D_1+2D_2+D_3
\end{eqnarray*}
For simplicity we assume $h>0$ in the following. The case $h<0$ can be
treated analogously.

We begin by estimating~$D_3$. Due to~\cite{roughI} one has
$\EX\|W_A(\cdot)\|^2 \in C^\infty(0,\infty)$. Hence, $|D_3|\leq C
\cdot h$ for some constant $C > 0$.

In order to estimate $D_2$, we define the shift-operator $\tau_t$ by
$\tau_t\om=\om(t+\cdot)$ for $t>0$.  Then the definitions of $\cG$ and
$\cF$ in~(\ref{Def-G}) and Lemma~\ref{f_global}, respectively, furnish
\begin{equation}     \label{Ft+h}
   \cG(\om)(t+h)
   = S(h)\cG(\om)(t)+\cF(\tau_t\om)(h) \; .
\end{equation}
Since $\{\om(s)\}_{s\in[0,t]}$ is stochastically independent of
$\{W(s)\}_{s\in[t,t+h]}$ we get
\begin{eqnarray*}
   \lefteqn{\EX < \cG(\om)(t+h),W_A(t+h)>}\\
   &=&\EX < S(h)\cG(\om)(t),\int_0^t S(t+h-s)d W(s)>\\
   & & \quad +
      \EX < \cF(\tau_t\om)(h),W_A(t+h)>\\
   &=&\EX < S(h)\cG(\om)(t),S(h)W_A(t)>+
      \EX < \cF(\tau_t\om)(h),W_A(t+h)>
\end{eqnarray*}
and together with the self-adjointness of $S(h)$ we finally arrive at
\begin{eqnarray}
  D_2 & = & \EX < (S(2h)-I)\cG(\om)(t),W_A(t)> \nonumber\\
  & & \quad + \; \EX < \cF(\tau_t\om)(h),W_A(t+h)> \; .
  \label{Gl-D2}
\end{eqnarray}
The boundedness of $\EX\|W_A(t)\|^2$ on~$J$ and Lemma~\ref{F-bound}
now yield
\begin{eqnarray*}
   |\EX < \cF(\tau_t\om)(h),W_A(t+h)>|
   &\le & \left(\EX \|\cF(\tau_t\om)(h)\|^2\right)^\frac12\cdot
     \left(\EX\|W_A(t+h)\|^2\right)^\frac12\\
   &\le& C \cdot h^{\frac14-\rho}.
\end{eqnarray*}
As for the first term in (\ref{Gl-D2}), notice that
\begin{equation}  \label{est-SG}
  \|(S(h)-I)v\| \leq
  \int_0^h \|(-A)S(s)v\| ds \leq
  C \cdot h^{a} \cdot \|(-A)^{a}v\|
\end{equation}
for any $v\in D((-A)^{a})$, with a constant $C$ which depends on $a
\in [0,1)$. Thus,
\begin{eqnarray*}
   \lefteqn{|\EX < (S(2h)-I)\cG(\om)(t),W_A(t)>|}\\
   &\le& C \cdot h^{\frac14-\rho} \cdot
   \left( \EX \left\| (-A)^{\frac14-\rho}\cF(\om)(t) \right\|^2
      + \EX \left\|(-A)^{\frac14-\rho}S(t)\om_0\right\|^2
      \right)^{\frac12}
\end{eqnarray*}
for any $\rho \in (0,1/4)$.  Together with
$\|(-A)^{\frac14-\rho}S(t)\om_0\| \le C \cdot \e^{\rho-\frac14} \cdot
\|\om_0\|$ and Lemma~\ref{F-bound} we eventually obtain
\begin{equation}\label{HoeErg1}
  |D_2|\le C \cdot h^{\frac14-\rho} \; .
\end{equation}
Finally we turn our attention to $D_1$.  Its definition and
(\ref{Ft+h}) imply
\begin{eqnarray*}
  \lefteqn{D_1 = \EX\|\cG(\om)(t+h)\|^2-\EX\|\cG(\om)(t)\|^2}\\
  &=& \EX < \cG(\om)(t+h)-\cG(\om)(t),
    \underbrace{\cG(\om)(t+h)+\cG(\om)(t)}_{=:{D_\cG}}>\\
  &=& \EX < (S(h)-I)\cG(\om)(t),D_\cG>
   +\EX <\cF(\tau_t\om)(h),D_\cG>.
\end{eqnarray*}
As in the discussion leading to (\ref{HoeErg1}), we obtain for any
$\tilde{a} \in [0,1/4)$
\begin{eqnarray}  \label{Hoe-G1}
  \lefteqn{|\EX < (S(h)-I)\cG(\om)(t),D_\cG>|}\nonumber\\
  &\leq& C \cdot h^{2\tilde{a}} \cdot
    \EX \left( \left\|(-A)^{\tilde{a}}\cG(\om)(t)\right\|
      \cdot \left\|(-A)^{\tilde{a}}D_\cG\right\| \right)
    \, \leq \, C \cdot h^{2\tilde{a}} \; .
\end{eqnarray}
Using again Lemma \ref{F-bound} we further derive
\begin{eqnarray}  \label{Hoe-G2}
  |\EX <\cF(\tau_t\om)(h),D_\cG>|
  &\leq& C \cdot \left( \EX\|\cF(\tau_t\om)(h)\|^2 \right)^\frac12
    \cdot \left( \EX\|\|D_\cG\|^2 \right)^\frac12\nonumber\\
  &\leq& C \cdot h^{\frac14-\rho} \; .
\end{eqnarray}
Combining (\ref{Hoe-G1}) for fixed $\tilde{a}$ near $1/4$ with
(\ref{Hoe-G2}) furnishes $|D_1| \leq C \cdot h^{\frac14-\rho}$, and we
finally obtain
\begin{displaymath}
  |Ens(t+h)-Ens(t)| \leq C \cdot h^{\frac14-\rho}
\end{displaymath}
for arbitrary $\rho \in (0,1/4)$.  This completes the proof of Theorem
\ref{holder}.

\section{Enstrophy Estimate: Asymptotics}

In Section~3 we established upper bounds on the growth of the
enstrophy $Ens(t) = \EX \|\om(t)\|^2 / 2$. Unfortunately, these bounds
fail to accurately describe the dynamics of $Ens(t)$ as $t \to 0$. For
example, the bound derived in Theorem~\ref{om-bound} will generally
not even converge to $Ens(0)$ as $t \to 0$. Therefore, this section is
devoted to investigating the small-time asymptotics of the enstrophy.
Similar to \cite{roughI, roughII} this will be accomplished by
relating $Ens(t)$ to the stochastic convolution.

In order to bound the growth of $\EX\|W_{A}(t)\|^2$, we assume that
the coefficients $\mu_k$ in (\ref{Wiener}) are bounded by $\mu_k^2 \le
c_\mu \cdot k^{-\delta}$ for some $\delta \in (0,1)$ and some positive
constant $c_\mu > 0$. In this situation we obtain similar to
\cite[Theorem~5.4]{roughI} the estimate
\begin{equation}  \label{linasym}
  \EX \|W_{A}(t)\|^2 \le C_0 \cdot t^\delta
\end{equation}
for arbitrary $t \in [0,T]$, where $C_0$ denotes a positive constant
which depends on~$T$.  Using the mild integral form (\ref{mild}) we
further get
\begin{displaymath}
  \|\om(t)-\om_0-W_A(t)\|
    \le \|(S(t)-I)\om_0\| + \|\cF(\om)(t)\| \; .
\end{displaymath}
If we now assume that $\EX\|\om_0\|^4<\infty$, then an application of
Lemma~\ref{F-bound} furnishes
\begin{eqnarray}  \label{asyG1}
\EX \|\om(t)-\om_0-W_A(t)\|^2
   &\le& 2\EX\|(S(t)-I)\om_0\|^2 + 2\EX\|\cF(\om)(t)\|^2\nonumber\\
   &\le& 2\EX\|(S(t)-I)\om_0\|^2 +C \cdot t^{\frac12-2\rho}\nonumber\\
   &\le& C \cdot t^{2\gamma} \cdot \EX\|(-A)^\gamma \om_0\|^2
      + C \cdot t^{\frac12-2\rho}
\end{eqnarray}
for fixed $\rho \in (0,1/4)$ and $\gamma \in [0,1)$. Thus, the
additional assumption $\EX\|(-A)^\gamma \om_0\|^2 < \infty$ for some
small $\gamma \in [0,1)$ implies
\begin{displaymath}
  \EX \|\om(t)-\om_0-W_A(t)\|^2 \le C \cdot t^{2\gamma} \; .
\end{displaymath}
Hence,
\begin{eqnarray*}
  \lefteqn{\left( \EX\|\om(t)\|^2 \right)^\frac12}\\
  &\le& \left( \EX\|\om_0\|^2 \right)^\frac12 +
    \left( \EX\|W_A(t)\|^2 \right)^\frac12 +
    \left( \EX\|\om(t)-\om_0-W_A(t)\|^2 \right)^\frac12\\
  & = & \left( \EX\|\om_0\|^2 \right)^\frac12 +
    {\cal O}\left(t^{\gamma}+t^{\frac{\delta}{2}}\right) \; .
\end{eqnarray*}
Similarly one obtains
\begin{eqnarray*}
  \lefteqn{\left( \EX\|\om(t)\|^2 \right)^\frac12}\\
  &\ge& \left( \EX\|\om_0\|^2 \right)^\frac12 -
    \left( \EX\|W_A(t)\|^2 \right)^\frac12 -
    \left( \EX\|\om(t)-\om_0-W_A(t)\|^2 \right)^\frac12\\
  & = & \left( \EX\|\om_0\|^2 \right)^\frac12 +
    {\cal O}\left( t^{\gamma} + t^{\frac{\delta}{2}} \right) \; ,
\end{eqnarray*}
and together these estimates show that $\EX\|\om(t)\|^2 =
\EX\|\om_0\|^2 + {\cal O}( t^{\gamma} + t^{\delta/2})$.  If on the
other hand we have $\om_0=0$, then (\ref{asyG1}) implies
\begin{displaymath}
  \EX \|\om(t)-W_A(t)\|^2 \le t^{\frac12-2\rho} \; ,
\end{displaymath}
which analogously results in $\EX\|\om(t)\|^2 = \EX\|W_A(t)\|^2 +
{\cal O}( t^{\frac14-\rho+\delta/2})$.  Using the definition of $Ens$,
this furnishes the following result on the small-time asymptotics of
the enstrophy.

\begin{theorem}[Asymptotics]
  Assume that $\EX\|(-A)^\gamma \om_0\|^2 < \infty$ for some small
  constant $\gamma>0$ and that $\EX\|\om_0\|^4<\infty$. Furthermore,
  suppose that (\ref{linasym}) holds for some small $\delta>0$. Then
  \begin{displaymath}
    Ens(t) = \frac12 \cdot \EX\|\om_0\|^2 +
      {\cal O}\left( t^{\gamma} + t^{\frac{\delta}{2}} \right)
    \; .
  \end{displaymath}
  If in addition we have $\om_0=0$ and let $\rho \in (0,1/4)$
  be arbitrary, then
  \begin{displaymath}
    Ens(t) = \frac12 \cdot \EX\|W_A(t)\|^2 +
      {\cal O}\left( t^{\frac{1+2\delta}{4}-\rho} \right) \; .
  \end{displaymath}
\end{theorem}

Notice that in the case $\om_0=0$ the second term on the right-hand
side is of higher order than $\EX\|W_A(t)\|^2 / 2$ only under
additional assumptions. For this we need $\delta < 1/2 - 2\rho$, as
well as a suitable lower bound on the growth of the first term
$\EX\|W_A(t)\|^2 / 2$ for small values of~$t$. The latter can be
achieved by imposing a lower bound on the growth of the
coefficients~$\mu_k$. For details we refer the reader to
\cite{roughI}.


\section{Summary}

The enstrophy $Ens(t)= \EX \|\om(t)\|^2 / 2$ is an averaged measure of
fluid vorticity $\om(t)$.  We have investigated the enstrophy evolution of
large-scale quasi-geostrophic flows under random wind forcing.
Thereby we have obtained results on upper bounds (Theorems 1 and 2),
H\"older continuity (Theorem 3), as well as small-time asymptotics
(Theorem 4) for the enstrophy.


\bigskip

{\bf Acknowledgements.}   
A part of this work was done at the
Oberwolfach Mathematical Research Institute, Germany,
while J. Duan was a Research
in Pairs Fellow, supported by {\em Volkswagen Stiftung}. This work
was partly supported by the NSF Grant DMS-9973204.


\end{document}